\documentclass[twoside,11pt]{article}

\usepackage{amsmath, amsfonts, amssymb, amsthm}
\usepackage{mathrsfs}
\usepackage{jmlr2e}

\firstpageno{1}

\begin{document}

\title{Interplay of minimax estimation and minimax support recovery under sparsity}

\author{\name Mohamed Ndaoud \email mohamed.ndaoud@ensae.fr \\
       \addr Department of Statistics\\
       CREST (UMR CNRS 9194), ENSAE\\
       5, av. Henry Le Chatelier, 91764 Palaiseau, FRANCE}

\editor{}

\maketitle

\begin{abstract}
In this paper, we study a new notion of scaled minimaxity for sparse estimation in high-dimensional linear regression model. We present more optimistic lower bounds than the one given by the classical minimax theory and hence improve on existing results. We recover sharp results for the global minimaxity as a consequence of our study. Fixing the scale of the signal-to-noise ratio, we prove that the estimation error can be much smaller than the global minimax error. We construct a new optimal estimator for the scaled minimax sparse estimation. An optimal adaptive procedure is also described.
\end{abstract}

 \section{Introduction}

\subsection{Statement of the problem}

Assume that we observe the vector of measurements $Y\in\mathbb{R}^{p}$ satisfying 
\begin{equation}\label{eq:def}
     Y = \beta + \sigma \xi 
\end{equation}
where $\beta \in\mathbb{R}^{p}$ is the unknown signal, $\sigma>0$ and the noise $\xi \sim \mathcal{N}\left(0,\mathbb{I}_{p} \right)$ is a standard Gaussian vector. Here, $\mathbb{I}_{p}$ denotes the $p\times p$ identity matrix. This model is a specific case of the more general model where $Y\in\mathbb{R}^{n}$ satisfies
\begin{equation}\label{eq:def2}
     Y = X\beta + \sigma \xi 
\end{equation}
where $X\in\mathbb{R}^{n \times p}$ is a given design or sensing matrix, and the noise is independent of $X$. Model \eqref{eq:def} corresponds to the orthogonal design. In this paper, we mostly focus on model \eqref{eq:def}. We denote by $\mathbf{P}_{\beta}$ the distribution of $Y$ in model \eqref{eq:def} or of $(Y,X)$ in model \eqref{eq:def2}, and by $\mathbf{E}_{\beta}$ the corresponding expectation.

We consider the problem of estimating the vector $\beta$. We will also explore its relation to the problem of recovering the support of $\beta$, that is the set $S_{\beta}$ of non-zero components of $\beta$. For an integer $s\le p$, we assume that $\beta$ is $s$-sparse, that is it has at most $s$ non-zero components.  We also assume that these components cannot be arbitrarily small. This motivates us to define the following set $\Omega^{p}_{s,a}$ of $s$-sparse vectors: 
$$  \Omega^{p}_{s,a} = \left\{\beta \in \mathbb{R}^{p}:\quad  |\beta|_{0} \leq s \quad \text{and} \quad   |\beta_{i}| \geq a, \quad \forall i \in S_{\beta} \right\}, $$
where $a>0$, $\beta_{i}$ are the components of $\beta$ for $i=1,\dots ,p,$ and $|\beta|_{0}$ denotes the number of non-zero components of $\beta$. The value $a$ characterizes the scale of the signal. In the rest of the paper, we will always denote by $\beta$ the vector to estimate, while $\hat{\beta}$ will denote the corresponding estimator. Let us denote by $\phi$ the scaled minimax risk
$$  \phi(s,a) =  \underset{\hat{\beta}}{\inf} \underset{\beta \in \Omega^{p}_{s,a}}{\sup}\mathbf{E}_{\beta}\left( \| \hat{\beta} - \beta\|^{2} \right), $$
where the infimum is taken over all possible estimators $\hat{\beta}$.
It is easy to check that $\phi$ is increasing with respect to $s$ and decreasing with respect to $a$. Note that, for $Y$ following model \eqref{eq:def}, the global minimax error over $\mathbb{R}^{p}$ is given by 
$$  \underset{\hat{\beta}}{\inf} \underset{\beta \in \mathbb{R}^{p}}{\sup}\mathbf{E}_{\beta}\left( \| \hat{\beta} - \beta\|^{2} \right)  = \sigma^{2} p. $$
The previous equality is achieved for $s=p$ and $a = 0$. Under the sparsity assumption, the previous result can be improved. In the seminal paper \cite{donoho92}, it is shown that the global sparse minimax estimation error has asymptotics:
\begin{equation}\label{eq:def3}
    \underset{\hat{\beta}}{\inf} \underset{ |\beta|_{0} \leq s}{\sup}\mathbf{E}_{\beta}\left( \| \hat{\beta} - \beta\|^{2} \right) = 2\sigma^{2} s \log(p/s)(1+o(1)) \quad \text{ as } \frac{s}{p}\to 0.
    \tag{$*$}
\end{equation}
Inspecting the proof of the minimax lower bound, one can see that \eqref{eq:def3} is achieved for $a = \sigma \sqrt{2\log{(p/s)}}(1+o(1))$.
We may also notice that the global sparse minimax estimation error is more optimistic than the global over $\mathbb{R}^{p}$. In this paper, we present an even more optimistic solution inspired by a notion of scaled minimax sparse estimation given by $\phi$. By doing so, we recover the global sparse estimation by taking the supremum over all $a$. In the rest of the paper, we will always denote by $SMSE$ the quantity $\phi$.

It is well known that minimax lower bounds are pessimistic. The worst case is usually specific to a critical region. Hence, a minimax optimal estimator can be good globally but may not be optimal outside of the critical region. By studying the quantity $\phi$ for fixed sparsity, we will emphasize this phenomenon. 

An optimistic lower bound for estimation of $s$-sparse vectors is given by $\sigma^{2}s$ and can be achieved when the support of vector $\beta$ is known. We say that an estimator $\hat{\beta}$ achieves exact estimation in model \eqref{eq:def} if
$$ \underset{s/p \to 0}{\lim}\underset{\beta \in \Omega_{s,a}^{p}}{\sup}\mathbf{E}_{\beta}\left( \| \hat{\beta} - \beta \|^{2} \right) = \sigma^{2}s. $$
We also say that estimator $\beta$ achieves exact support recovery in model \eqref{eq:def} if
$$ \underset{p \to \infty }{\lim} \underset{\beta \in \Omega_{s,a}^{p}}{\sup}\mathbf{P}_{\beta}\left( S_{\hat{\beta}} = S_{\beta} \right) = 1. $$
In this paper, we shed some light on the relation between exact support recovery and exact estimation. Specifically, we give an answer to the following questions that motivate the present work.
\begin{itemize}
    \item How pessimistic is the result \eqref{eq:def3}? Can we do any better by fixing the scale value $a$?
    \item Is exact support recovery necessary to achieve exact estimation?
    \item Can we achieve minimax optimality with respect to SMSE adaptively to the scale value $a$?
\end{itemize}
In the dense regime where $s \asymp p$, the minimax estimation error is of order $\sigma^{2}p$ independently of $a$. Hence, in the rest of the paper, we focus on the regime where $\frac{s}{p} = o(1)$. All the proofs are deferred to Appendix.

{\bf Notation.}
In the rest of this paper we use the following notation. For given sequences $a_{n}$ and $b_n$, we say that $a_{n} = \mathcal{O}(b_{n})$ (resp $a_{n} = \Omega(b_n)$) when $a_{n} \leq c b_{n}$ (resp $a_{n} \geq c b_{n}$) for some absolute constant $c>0$. We write $a_{n} \asymp b_{n} $ if $a_{n} = \mathcal{O}(b_{n})$ and $a_{n}=\Omega(b_{n})$. For ${\bf x},{\bf y}\in \mathbb{R}^{p}$, $\|{\bf x}\|$ is the Euclidean norm of ${\bf x}$, and ${\bf x}^{\top}{\bf y}$ the corresponding inner product. For $q \geq 1$, and ${\bf x}\in \mathbb{R}^{p}$, we denote by $\| {\bf x}\|_{q}$ the $l_{q}$ norm of ${\bf x}$. For a matrix $X$, we denote by $X_j$ its $j$th column. For $x,y \in \mathbb{R}$, we denote by $x\vee y$ the maximum of $x$ and $y$ and we set $x_{+}=x \vee 0$. For $q \geq 1$ and $\xi$ a centered Gaussian random variable with variance $\sigma^{2}$, we denote by $\sigma_{q}$ the quantity $\mathbf{E}(|\xi|^{q})^{1/q}$.
The notation  $\mathbf{1}(\cdot)$ stands for the indicator function. We denote by $C$ and $C_{q}$ positive constants where the second one depends on $q$ for some $q \geq 1$.
\subsection{Related literature}
The literature on minimax sparse estimation in high-dimensional linear regression (for both random and orthogonal design) is very rich and its complete overview falls beyond the format of this paper. We mention here only some recent results close to our work. All sharp results are considered in the regime $\frac{s}{p} \to 0$.
\begin{itemize}
    \item In \cite{bellec}, the authors show that SLOPE estimator, which is defined in \cite{candes}, is minimax optimal for estimation under sparsity constraint in model \eqref{eq:def2}, as long as $X$ satisfies some general conditions. This result is non-asymptotic.
    \item \cite{bellec2018} proves that the minimax estimation rate of convex penalized destimators cannot be improved for sparse vectors, even when the scale parameter $a$ is large. This fact is mainly due to the bias caused by convex penalization as it is the case for LASSO and SLOPE estimators.
    \item In \cite{su2016}, it is shown that SLOPE is sharply minimax optimal on $\{ |\beta |_{0} \leq s\}$ giving the asymptotic optimal estimation error $2\sigma^{2}s\log{\frac{p}{s}}$ in model \eqref{eq:def}. In model \eqref{eq:def2}, where $X$ has i.i.d standard Gaussian entries and under the asymptotic condition $\frac{s\log{p}}{n} \to 0$, SLOPE gives the asymptotic optimal error $2\frac{\sigma^{2}}{n}s\log{\frac{p}{s}}$, cf. \cite{su2016}.
    \item \cite{zhou} show that the penalized least squares estimator with a penalty that grows like $2\sigma^{2}s\log{\frac{p}{s}}$, is sharply minimax optimal on $\{ |\beta|_{0} \leq s \}$ under additional assumptions on $s$ and $p$.
\end{itemize}
\subsection{Main contribution}
Inspired by the related literature, the present work is also motivated by the following questions.
\begin{itemize}
    \item In model \eqref{eq:def}, the proof of lower bounds uses a worst case vector with non-zero components that scale as $\sigma \sqrt{2\log{\frac{p}{s}}}$ in order to get the best lower bound. In other words, the worst case happens for a specific vector $\beta$. Can we do better far from this vector?
    \item One of the popular approaches is to recover the support of a sparse vector and then estimate this vector on the obtained support. In this case the error of estimation is of order $s \sigma^{2}$ and is the best one can hope to achieve. Is it necessary to recover the true support in order to get this error? This is an important question that we address in this paper.
    \item If the answer to the previous question is negative, can we propose an algorithm that would be optimal in the sense of SMSE, practical and adaptive?
\end{itemize}

The main contribution in this paper is a sharp study of the minimax risk $\phi$. What is more, we study a more general quantity given by
\begin{equation}\label{eq:minimax}
    \underset{\hat{\beta}}{\inf} \underset{\beta \in \Omega^{p}_{s,a}}{\sup}\mathbf{E}_{\beta}\left( \| \hat{\beta} - \beta\|_{q}^{q} \right),
\end{equation}
for any $q \geq 1$. We give lower bounds and corresponding upper bounds for \eqref{eq:minimax}. We show that in the interesting regime $\frac{s}{p} = o(1)$, our lower and upper bounds match not only in the rate but also in the constant up to a factor 4 under a mild condition on sparsity. As a result of our study, we recover two interesting phase transitions when $\frac{s}{p} = o(1)$. 

The first one is that there are basically two regimes in estimation. 
For $a \leq \sqrt{2\sigma^{2}\log(p/s)}(1-\epsilon) $ and $\epsilon>0$, the asymptotic SMSE is $2s\sigma^{2}\log(p/s)$. This regime is called the hard recovery regime, where we prove that the error is due to misspecification in recovering the support. Alternatively, for $a \geq (1+\epsilon)\sqrt{2\sigma^{2}\log(p/s)}$, the error is of order $s\sigma^{2}$. This regime is presented as the hard estimation regime. In this regime, we can recover a good fraction of the support but still have to pay for the estimation on the support. Hence, and surprisingly, the SMSE is almost piece-wise constant as a function of $a$. This shows that the sparse minimax risk can be made much smaller once we get far from some critical region. 

Another contribution of this paper is a new phase transition related to sparsity. In \cite{butucea2018}, it is shown that a necessary condition to achieve exact recovery is given by $a \geq \sigma\sqrt{2\log(p-s)} + \sigma\sqrt{2\log(s)}$.To achieve exact estimation, a necessary condition is $a \geq \sigma\sqrt{2\log(p/s-1) + 2\log\log(p/s-1)} + \sigma\sqrt{2\log\log(p/s-1)} $. Hence exact recovery is not necessary for exact estimation. In fact, when $s \gg \log(p)$ then exact estimation is easier and when $s \ll \log(p)$ exact recovery becomes easier. This shows that there is no direct implication of exact recovery on exact estimation, hence the latter task should be considered as a separate problem. 

Finally, one more contribution of this paper is adaptivity. We give an optimal adaptive variant of our procedure, that achieves the sparse minimax optimal rate and whenever exact estimation is possible achieves it as well. By doing so, our procedure improves on the existing literature. In fact, Lasso is known to have an unavoidable bias of order $\sigma^{2}s\log(p/s)$ even on the class $\Omega_{s,a}^{p}$, cf. \cite{bellec2018}. We show that our procedure is better in the sense that it gets rid of the bias whenever it is possible.
\section{Towards more optimistic lower bounds for estimation}\label{sec:nas}
In several papers, lower bounds for minimax risk are derived using the Fano lemma. These lower bounds are usually far from being sharp in the non-asymptotic setting. We establish, in this section, non-asymptotic lower bounds on the minimax risk based on some revisited two-hypothesis testing techniques. 

We derive two lower bounds for the SMSE. The scaled error of estimation of sparse vectors can be decomposed into two types of error. A first one based on the error of estimation when the true support $S_{\beta}$ is known and a second one is given by the error of recovery of the true support when the vector components are known but not the support. For this purpose, we prove first a general lower bound for constrained minimax sparse estimation.

In the next theorem, we reduce the constrained minimax risk over all estimators to a Bayes risk with arbitrary prior measure $\pi$ on $\mathbb{R}^p$ and give a bound on the difference between the two risks. This result is true in a general setup, non necessarily for Gaussian models. For a particular choice of measure $\pi$,  we provide an explicit bound of the remainder term.

Consider the set of vectors $\Theta_{s,a} \subseteq \mathbb{R}^{p}$, and assume that we are given a family $\{P_{\beta}, \beta \in \Theta_{s,a}\}$ where each $P_{\beta}$ is a probability distribution on a measurable space $(\mathcal{X}, \mathcal{U})$.
We observe $Y$  drawn from $P_{\beta}$ with some unknown $\beta \in \Theta_{s,a}$ and we consider the risk of an estimator $\hat{\beta}=\hat{\beta}(Y)$:
$$
\sup_{ \beta \in \Theta_{s,a}}\mathbf{E}_{\beta}\|\hat{\beta} - \beta\|_{q}^{q}
$$
where $\mathbf{E}_{\beta}$ is the expectation with respect to $P_{\beta}$. Let $\pi$ be a probability measure on $\mathbb{R}^{p}$ (a prior on $\beta$). We denote by $\mathbb{E}_{\pi}$ the expectation with respect to~$\pi$.
\begin{theorem}\label{th:low1}
For any $s<p$, $q\geq1$ and any probability measure $\pi$ on $\mathbb{R}^{p}$, there exists $C_{q}>0$ such that
\begin{equation}\label{eeq1:th}
\inf_{\hat{\beta}} \sup_{ \beta \in \Theta_{s,a}} \mathbf{E}_\beta \|\hat{\beta} - \beta \|_{q}^{q}
\geq  \inf_{\hat T\in \mathbb{R}^p}\mathbb{E}_{\pi}\mathbf{E}_\beta \sum_{j=1}^p | {\hat T}_j(Y) - \beta_j |^{q} - C_{q}\,\mathbb{E}_{\pi}\left[ \left( \mathbf{E}(\|\beta^{A}\|^{q}_{q}|Y) + \|\beta\|_{q}^{q} \right)\mathbf{1}( \beta \not\in \Theta_{s,a} ) \right] ,
\end{equation}
where $\beta^{A}:= \beta.\mathbf{1}(\beta \in \Theta_{s,a})=(\beta_{1} \mathbf{1}( \beta\in \Theta_{s,a} ),\dots,\beta_{p} \mathbf{1}( \beta\in \Theta_{s,a} ) )$, $\inf_{\hat{\beta}}$ is the infimum over all estimators and $\inf_{\hat T\in \mathbb{R}^p}$ is the infimum over all  estimators $\hat T(Y)=({\hat T}_1(Y),\dots,{\hat T}_p(Y))$ with values in $\mathbb{R}^p$.
\end{theorem}
Theorem \ref{th:low1} is valid in a very general setting. We present now specific lower bounds in the general model of linear regression. Assume that $Y \in \mathbb{R}^{n}$ follows model \eqref{eq:def2}, where $X$ is a deterministic design. The following lemma is useful to get more precise lower bounds in model \eqref{eq:def2}. It is based on the simple observation that under independent prior distributions of the entries of $\beta$ the oracle estimator of a given component does not depend on the rest of the components.
\begin{lemma}\label{lem:low}
Assume that $Y$ satisfies model \eqref{eq:def2} with a deterministic design $X$. Then 
$$  \inf_{\hat T\in \mathbb{R}^p}\mathbb{E}_{\pi}\mathbf{E}_\beta \sum_{j=1}^p | {\hat T}_j(X,Y) - \beta_j |^{q} \geq \sum_{j=1}^{p}\inf_{\hat T_{j}\in \mathbb{R}}\mathbb{E}_{\pi_{j}}\mathbf{E}_{\beta_{j}} | {\hat T}_j(X_{j},\tilde{Y}_{j}) - \beta_j |^{q} , $$
where $\tilde{Y}_{j} = Y - \sum_{i\neq j} X_{i}\beta_{i} = \beta_{j}X_{j} + \sigma \xi$.
\end{lemma}
Using the previous lemma, we are now ready to give two sharp lower bounds for the SMSE. A first one supposed to capture the error of estimation when the support is known, while the second one handles the case where the support is not known.
\begin{theorem}\label{th:low2}
Assume that $Y$ follows model \eqref{eq:def2} with a deterministic design $X$. For any $a>0$, $q \geq 1$ and $s<p$ we have 
$$ \underset{\hat{\beta}}{\inf}\underset{\beta \in \Omega^{p}_{s,a}}{\sup} \mathbf{E}_{\beta}\left( \| \hat{\beta} - \beta \|_{q}^{q}\right)\geq \sigma_{q}^{q} \underset{|S|= s}{\max}\sum_{i \in S}\frac{1}{\|X_{i}\|^{q}_{2}}.$$
\end{theorem}
In order to derive the next lower bound, we define the quantity $\Psi$ introduced in \cite{ndaoud2018} in the context of support recovery:
$$ \Psi(p,s,a,\sigma,X) := \sum_{j=1}^{p}\left(\frac{s}{p}\mathbf{P}(\sigma \varepsilon \geq (a-t_{j}(a))\|X_{j}\|_{2})  + (1-\frac{s}{p})\mathbf{P}(\sigma \varepsilon \geq t_{j}(a)\|X_{j}\|_{2})  \right), $$
where $\varepsilon$ is standard Gaussian random variable and
$$ t_{j}\left( a \right) := \frac{a }{2} + \frac{\sigma^{2}\log\left(\frac{p}{s}-1\right)}{a\|X_{j}\|_{2}^{2}}, \quad \forall j=1,\dots ,p.$$
\begin{theorem}\label{th:low3}
Assume that $Y$ follows model \eqref{eq:def2} with deterministic design $X$. For any $a>0$, $q\geq1$ and $s<p$ we have 
$$ \forall s' \in (0,s), \quad \underset{\hat{\beta}}{\inf}\underset{\beta \in \Omega^{p}_{s,a}}{\sup} \mathbf{E}_{\beta}\left( \| \hat{\beta} - \beta \|_{q}^{q}\right)\geq a^{q}\frac{s'}{s}\left(\frac{1}{2^{q}}\Psi(p,s,a,\sigma,X) - 2se^{-\frac{(s-s')^{2}}{2s}}\right).$$
\end{theorem}

The proof is based on arguments similar to \cite{butucea2018}.
Assume now that we are under model \eqref{eq:def} and set 
$$  \psi\left(p,s,a,\sigma \right) := \left(p-s\right)\mathbf{P}\left( \sigma \varepsilon > t\left(a\right) \right) + s \mathbf{P}\left( \sigma \varepsilon > a - t(a)\right),$$
where $\varepsilon$ is a standard Gaussian random variable and
\begin{equation}\label{eq:thres1}
    t\left( a \right) := \frac{a }{2} + \frac{\sigma^{2}\log\left(\frac{p}{s}-1\right)}{a}. 
\end{equation}
The minimax Hamming loss for model \eqref{eq:def} was studied in \cite{butucea2018}, where it was shown that it is very linked to $\psi$. One may notice that, under model \eqref{eq:def}, $\Psi(p,s,a,\sigma,\mathbb{I}_{p}) = \psi(p,s,a,\sigma)$. We define now the following estimation rate
$$
\Phi(a) := \left\{
    \begin{array}{ll}
        a^{q}\psi(s,p,a,\sigma) \vee \sigma_{q}^{q}s & \mbox{if } a \geq t^{*},\\
        s\sigma^{q} \left(2\log(\frac{p}{s}-1)\right)^{\frac{q}{2}} & \mbox{else,}
    \end{array}
\right.
$$
where 
\begin{equation}\label{eq:thres2}
    t^{*}=\sigma \sqrt{2\log{\frac{p}{s}-1}}. 
\end{equation}
The next proposition is a consequence of previous theorems and shows the link between the minimax Hamming loss and the minimax estimation risk.
\begin{proposition}\label{prop:low}
Assume that $Y$ follows model \eqref{eq:def}. For any $a>0$, $q\geq1$, $s<p/2$ and $s \geq 8q\log\log(p)$, we have 
$$ \underset{\hat{\beta}}{\inf}\underset{\beta \in \Omega^{p}_{s,a}}{\sup} \mathbf{E}_{\beta}\left(\|\hat{\beta} - \beta\|^{q}_{q} \right)\geq C_{q}\Phi(a),$$
where $C_{q}>0$. 
\end{proposition}
\begin{remark}
The mild condition $s = \Omega( \log\log{p})$ is an artifact of the proof of the lower bound. We believe that this condition can be removed or further relaxed.
\end{remark}
A more careful proof of the previous result can lead us to $C_{q}=(1+o(1))$ as $\frac{s}{p} \to 0$. We omit the proof of this, since we give a more accurate result in the next section. Analyzing the lower bound of Proposition \ref{prop:low}, it turns out that the minimax rate $\sigma^{2}s\log{(p/s)}$, for $q=2$, cannot be improved when $a \leq t^{*}$. We will see later that this is not the case for large values of $a$. The next section is devoted to closing this gap by deriving matching upper bounds.
\section{Optimal scaled minimax estimators}\label{sec:3}
In this section, we consider upper bounds for the scaled minimax risk under model \eqref{eq:def}. For $a>0$ define the following estimator:
\begin{equation}\label{eq:est1}
    \hat{\beta}^{a}_{j} := Y_{j}\mathbf{1}_{\{|Y_{j}| \geq t(a \vee t^{*})\}},  \quad \forall j \in 1,\dots,p,
\end{equation}
where $t(.)$ and $t^{*}$ are defined respectively in \eqref{eq:thres1} and \eqref{eq:thres2}. The following result gives a matching upper bound for the scaled minimax risk.
Set
$$
\Phi_{+}(a) := \left\{
    \begin{array}{ll}
        a^{q}\psi_{+}(p,s,a,\sigma) \vee \sigma_{q}^{q}s & \mbox{if } a \geq t^{*},\\
        s\sigma^{q} \left(2\log(\frac{p}{s}-1)\right)^{\frac{q}{2}} & \mbox{else,}
    \end{array}
\right.
$$
where $\psi_{+}$ is given by
$$ \psi_{+}\left(p,s,a,\sigma \right) := \left(p-s\right)\mathbf{P}\left( \sigma \varepsilon > t\left(a\right) \right) + s \mathbf{P}\left( \sigma \varepsilon > (a - t(a)_{+})\right),$$
and $\varepsilon$ is a standard Gaussian random variable. Notice that $\Phi_{+}(a) \leq \Phi(a)$. This remark, combined with the next theorem, shows minimax optimality of the estimator \eqref{eq:est1}.
\begin{theorem}\label{thm:up1}
Assume that $Y$ follows model \eqref{eq:def}. For all $a>0$, let $\hat{\beta}^{a}$ be the estimator \eqref{eq:est1}. For all $q \geq 1$ and $s< p/2$ we have
$$  \underset{\beta \in \Omega^{p}_{s,a}}{\sup} \mathbf{E}_{\beta}\left( \| \hat{\beta}^{a} - \beta \|_{q}^{q}\right)\leq C_{q}\Phi_{+}(a),$$
where $C_{q}$ is a universal constant depending only in $q$.
\end{theorem}
Combining this result with Proposition \ref{prop:low}, we deduce the next corollary.
\begin{corollary}\label{cor:up1}
Assume that $Y$ follows model \eqref{eq:def}. For all $a>0$, let $\hat{\beta}^{a}$ be the estimator \eqref{eq:est1}. For all $q \geq 1$, $s<p/2$ and $s \geq 8q\log\log(p)$, there exists $C_{q}>0$ such that
$$ \underset{\beta \in \Omega^{p}_{s,a}}{\sup} \mathbf{E}_{\beta}\left( \| \hat{\beta}^{a} - \beta \|_{q}^{q}\right)\leq C_{q}\underset{\hat{\beta}}{\inf}\underset{\beta \in \Omega^{p}_{s,a}}{\sup} \mathbf{E}_{\beta}\left( \| \hat{\beta} - \beta \|_{q}^{q}\right).$$
\end{corollary}
We give now a more accurate upper bound in the regime $\frac{s}{p} \to 0$. Assume that $s \leq p/4$. For $ q\geq 1$, and $\epsilon \in [0,1]$ define 
$$ a_{q}(\epsilon) = \sqrt{2\sigma^{2}\log(\frac{p}{s}-1) + q \epsilon \sigma^{2}\log\log(\frac{p}{s}-1)} + \sqrt{ q \epsilon \sigma^{2}\log\log(\frac{p}{s}-1)}. $$
Set 
$$
\Phi_{o}(a) := \left\{
    \begin{array}{ll}
        s\sigma^{q} \left(2 \log(\frac{p}{s}-1)\right)^{\frac{q}{2}}  & \mbox{if } a \leq a_{q}(0),\\
        \frac{s\sigma^{q} \left(2 \log(\frac{p}{s}-1)\right)^{\frac{q}{2}(1-\epsilon)}}{1 + \sigma\sqrt{\frac{\pi}{2}\epsilon q \log\log(\frac{p}{s}-1)}} \vee \sigma_{q}^{q}s & \mbox{if } a = a_{q}(\epsilon) , \epsilon \in (0,1),\\
        s\sigma_{q}^{q} & \mbox{if } a \geq a_{q}(1).\\
        
    \end{array}
\right.
$$
The next theorem gives sharp upper bounds in the regime $\frac{s}{p}\to 0$.
\begin{theorem}\label{thm:up2}
Assume that $Y$ follows model \eqref{eq:def}. For all $a>0$, let $\hat{\beta}^{a}$ be the estimator \eqref{eq:est1}. In the regime where $\frac{s}{p} \to 0$, for all $q \geq 1$, we have
$$
\underset{\beta \in \Omega^{p}_{s,a}}{\sup} \mathbf{E}_{\beta}\left( \| \hat{\beta}^{a} - \beta \|_{q}^{q}\right) \leq \Phi_{o}(a)(1+o(1)).
$$
\end{theorem}
As a consequence of previous results, we derive the next corollary that gives an almost sharp bound for SMSE when $\frac{s}{p} \to 0$.
\begin{corollary}\label{cor:last}
Assume that $Y$ follows model \eqref{eq:def}. For all $a>0$, $q \geq 1$ and $s \geq 8q\log\log(p)$, in the regime $\frac{s}{p} \to 0$, we have
$$\frac{1}{4} +o(1) \leq \underset{\hat{\beta}}{\inf}\underset{\beta \in \Omega^{p}_{s,a}}{\sup} \frac{ \mathbf{E}_{\beta}\left( \| \hat{\beta} - \beta \|_{q}^{q}\right)}{\Phi_{o}(a)} \leq 1+o(1),$$
and
$$
\underset{\hat{\beta}}{\inf}\underset{\beta \in \Omega^{p}_{s,a}}{\sup} \frac{ \mathbf{E}_{\beta}\left( \| \hat{\beta} - \beta \|_{q}^{q}\right)}{s\sigma_{q}^{q}} = 1+o(1) \quad \text{if } a \geq a_{q}(1).
$$
\end{corollary}
Inspecting the proof of Corollary \ref{cor:last}, we may notice that the discrepancy between the bounding constants is mainly caused by values of the scale $a=a_{q}(\epsilon)$ such that $\epsilon \in (0,1)$. Corollary \ref{cor:last} shows that we can construct an almost sharp optimal minimax estimator provided $a$ and $s$. The next section is devoted to the question of adaptivity.

\section{ Adaptative scaled minimax estimators }

In Section \ref{sec:3} we have shown that the minimax rate is given by the quantity $\Phi_{o}(a)$ in a sharp way if $\frac{s}{p} \to 0$. Note that $\Phi_{o}(a)$ is almost piece-wise constant as a function of $a$. In fact the study of $\Phi_{o}(a)$ gives rise to three different regimes that we describe below.
\begin{enumerate}
    \item \textit{Hard recovery regime:}
    
Let $a \leq \sqrt{2\sigma^{2}\log\left(\frac{p}{s}-1\right)}$. We call this the hard recovery regime. In this regime, $\Phi_{o}(a)$ is constant and has a value of order $ \sigma^{q}s\left(2\log\left(\frac{p}{s}-1\right)\right)^{q/2}$. It turns out that the worst case of estimation happens for $a = \sqrt{2\sigma^{2}\log\left(\frac{p}{s}-1\right)}$. This error is mainly due to the fact the we cannot achieve almost full recovery as defined in \cite{butucea2018}.
    \item \textit{Hard estimation regime:}
    
This regime corresponds to values of $a$ such that $a \geq \sigma\sqrt{2\log\left( \frac{p}{s}-1\right)} \sqrt{1+4\frac{q\log\log(\frac{p}{s}-1)}{\log(\frac{p}{s}-1)}} .$ In this regime $\Phi_{o}(a)$ is of order $\sigma_{q}^{q} s$. In this region, the error of estimation on a known support dominates the error of recovering the support.

\item \textit{Transition regime:}

This regime concerns the remaining values of $a$ falling between the two previous regimes. In this regime $\Phi_{o}(a)$ is not constant any more. It represents a monotonous and continuous transition from one regime to another. 
\end{enumerate}
After analyzing the SMSE, we give a couple of remarks.
\begin{remark}

\begin{itemize}

    \item If $s = \small{o}(p) $ there are basically two regimes around the threshold $\sqrt{2\sigma^{2}\log\left( \frac{p}{s}-1\right)}$. Notice also that the hard estimation error is very small compared to the hard recovery error. We may notice that the SMSE is very small compared to the minimax sparse estimation error in the hard estimation regime. This proves how pessimistic the general minimax lower bounds are and that we can do much better for the scaled minimax risk.
    \item The case $s \sim p$ is of small interest. There is no phase transition in this case, since the SMSE is of order $\sigma^{q} p$ for every $a$.
    \item In the Hard estimation regime, the minimax error rate is the same as if the support were exactly known. It is interesting to notice that we need a weaker condition to get this rate when $s \gg \log(p) $, while a stronger necessary condition is needed for exact recovery, cf. \cite{butucea2018}. Hence exact support recovery is not necessary to achieve exact estimation.
\end{itemize}
\end{remark}
Notice also that the transition regime happens in a very small neighborhood around the universal threshold $\sqrt{2\sigma^{2}\log{\frac{p}{s}}}$.
Thus, it is very difficult to be adaptive to $a$ in the transition regime. For $s \leq p/4$, define the following estimator:
\begin{equation}\label{eq:est2}
    \hat{\beta}^{s}_{j} := Y_{j}\mathbf{1}_{\{|Y_{j}| \geq t^{*}_{s}\}},  \quad \forall j \in 1,\dots,p, 
\end{equation}
where 
$$t^{*}_{s} := \sqrt{2\sigma^{2}\log(p/s-1) + \sigma^{2}q\log\log(p/s-1)}.$$
We define a more convenient adaptive estimation error.
Set
$$
\Phi_{ad}(a) := \left\{
    \begin{array}{ll}
         \sigma_{q}^{q}s & \mbox{if } a \geq a_{q}(1),\\
        s\sigma^{q} \left(2\log(\frac{p}{s}-1)\right)^{\frac{q}{2}} & \mbox{else.}
    \end{array}
\right.
$$

The following result gives a matching upper bound for the adaptive scaled minimax risk.
\begin{theorem}\label{thm:up1ad}
Assume that $Y$ follows model \eqref{eq:def}. Let $\hat{\beta}^{s}$ be the estimator \eqref{eq:est2}. For all $q \geq 1$, $a>0$ and $s< p/4$ we have
$$  \underset{\beta \in \Omega^{p}_{s,a}}{\sup} \mathbf{E}_{\beta}\left( \| \hat{\beta}^{s} - \beta \|_{q}^{q}\right)\leq C_{q}\Psi_{ad}(a),$$
where $C_{q}$ is a universal constant depending only in $q$.
\end{theorem}
Since the two main regimes are hard estimation and hard recovery, we restricted the notion of adaptivity to these regimes. By doing so, we constructed an almost optimal estimator adaptively to the parameter $a$. This estimator is minimax optimal over the set of $s$-sparse vectors and achieves exact estimation when necessary conditions are satisfied. Our estimator has a phase transition around the universal threshold. Based on a procedure similar to \cite{butucea2018}, we can also construct an optimal estimator adaptive to sparsity. We do not give further details here for the sake of brevity. 
\section{Conclusion:}
In this paper, we define and study a new notion that we call scaled minimax sparse estimation. We assess how pessimistic are minimax lower bounds for the problem of sparse estimation. We also show that exact recovery is not necessary for exact estimation in general. As a result, we construct a new estimator optimal for the SMSE and present its adaptive version, improving on existing procedures for the problem of estimation. 

\acks{We would like to thank Alexandre Tsybakov for valuable comments on early
versions of this manuscript.}


\newpage

\appendix
\section*{Appendix}
\label{app:theorem}
The following bounds for the tails of Gaussian distribution will be useful:
$$ \frac{e^{-y^{2}/2}}{\sqrt{2\pi}y + 4 } \leq \frac{1}{\sqrt{2\pi}}\int_{y}^{\infty}e^{-u^{2}/2}du \leq \frac{e^{-y^{2}/2}}{\sqrt{2\pi}y \vee 2}.$$
for all $y \geq 0$. These bounds are an immediate consequence of formula 7.1.13 in
\cite{Gauss} with $x = y/\sqrt{2}$.
\begin{proof}[\textbf{Proof of Theorem \ref{th:low1}}]
Throughout the proof, we write for brevity $A=\Theta_{s,a}$.
Set $ \beta^{A} = \beta. \mathbf{1}(\beta \in A)$ and denote by $\pi_A$ the probability measure $\pi$ conditioned by the event $\{\beta \in A\}$,  that is, for any $C\subseteq \mathbb{R}^{d}$,
$$
\pi_A(C) = \frac{\pi (C\cap \{\beta \in A\})}{\pi(\beta \in A)}\,.
$$
The measure $\pi_A$ is supported on $A$ and we have
\begin{eqnarray*}
\inf_{\hat{\beta}} \sup_{ \beta \in A} \mathbf{E}_\beta |\hat{\beta} - \beta|^{q}_{q}
&\ge&
 \inf_{\hat{\beta}} \mathbb{E}_{\pi_A}\mathbf{E}_\beta |\hat{\beta} - \beta|^{q}_{q}
= \inf_{\hat{\beta}} \mathbb{E}_{\pi_{A}}\mathbf{E}_\beta |\hat{\beta} - \beta^A|^{q}_{q}
\\
&\ge&   \sum_{j=1}^{p} \inf_{\hat{T_j}} \mathbb{E}_{\pi_A}\mathbf{E}_\beta |\hat{T}_j - \beta^A_j |^{q}
\end{eqnarray*}
where $\inf_{\hat T_j}$ is the infimum over all estimators $\hat T_j= \hat T_j(Y)$ with values in $\mathbb{R}$.
According to
Theorem~1.1 and Corollary~1.2 on page 228 in \cite{LC}, there exists a Bayes estimator $B^A_j=B^A_j(Y)$ such that
$$
\inf_{\hat{T_j}} \mathbb{E}_{\pi_A}\mathbf{E}_\beta |\hat{T}_j - \beta^A_j |^{q}=\mathbb{E}_{\pi_A}\mathbf{E}_\beta |B^A_j - \beta^A_j |^{q}.
$$
In particular, for any estimator $\hat T_j(Y)$ we have
\begin{equation}\label{eeq0}
\mathbb{E}^A\big(|B^A_j(Y) - \beta^A_j |^{q} \big| Y\big) \le \mathbb{E}^A\big(|\hat T_j(Y) - \beta^A_j |^{q} \big| Y\big)
\end{equation}
almost surely. Here, the superscript $A$ indicates that the conditional expectation $\mathbb{E}^A(\cdot|Y)$ is taken when $\beta$ is distributed according to $\pi_A$.
Therefore,
\begin{equation}\label{eeq1}
\inf_{\hat{\beta}} \sup_{ \beta \in A} \mathbf{E}_\beta |\hat{\beta} - \beta|_{q}^{q}
\geq \mathbb{E}_{\pi_A}\mathbf{E}_\beta \sum_{j=1}^{p} | B^A_j - \beta^{A}_{j}|^{q} .
\end{equation}
Using this, we obtain
\begin{align}\nonumber
  \inf_{ \hat{T}\in\mathbb{R}^p} \mathbb{E}_{\pi}\mathbf{E}_\beta | \hat{T} - \beta |^{q}_{q} & \leq \mathbb{E}_{\pi}\mathbf{E}_\beta \sum_{j=1}^{p} | B^A_j - \beta_{j}|^{q} \\ \nonumber
    & = \mathbb{E}_{\pi}\mathbf{E}_\beta  \Big(\sum_{j=1}^{p} |B^A_j - \beta_{j}|^{q} \mathbf{1}(\beta \in A) \Big) + \mathbb{E}_{\pi}\mathbf{E}_\beta
   \Big(\sum_{j=1}^{p} |B^A_j - \beta_{j}|^{q}  \mathbf{1}(\beta \not\in A) \Big)
  \\  \nonumber
  & \leq \mathbb{E}_{\pi_A}\mathbf{E}_\beta  \sum_{j=1}^{p} |B^A_j - \beta_{j}^A|^{q}  + \mathbb{E}_{\pi}\mathbf{E}_\beta
   \Big(\sum_{j=1}^{p} |B^A_j - \beta_{j}|^{q}  \mathbf{1}(\beta \not\in A)  \Big)
   \\  \label{eeq2}
  & \le \mathbb{E}_{\pi_A}\mathbf{E}_\beta  \sum_{j=1}^{p} |B^A_j - \beta_{j}^A|^{q}  +
  \mathbb{E}_{\pi}\mathbf{E}_\beta
   \sum_{j=1}^{p} 2^{q-1}(|B^A_j|^{q}  + |\beta_j|^{q}) \mathbf{1}(\beta \not\in A) .
 \end{align}
Our next step is to bound the term
$$
\mathbb{E}_{\pi}\mathbf{E}_\beta
   \sum_{j=1}^{p} |B^A_j|^{q} \mathbf{1}(\beta \not\in A).
$$
For this purpose, we first note that inequality \eqref{eeq0} with $\hat T_j(Y) = 0$ implies that
$$
|B^A_j (Y)|^{q} = \mathbb{E}^A(|B_j^A(Y)|^{q}|Y) \le  2^{q}\mathbb{E}^A(|\beta_j^A|^{q}|Y).
$$
Thus
$$
\mathbb{E}_{\pi}\mathbf{E}_\beta
    \sum_{j=1}^{p} |B^A_j|^{q} \mathbf{1}(\beta \not\in A)\le
   2^{q}\mathbb{E}_{\pi}
     \mathbb{E}^A(\|\beta^A\|_{q}^{q}|Y)\mathbf{1}(\beta \not\in A).
$$
Combining this inequality with \eqref{eeq1} and \eqref{eeq2} yields \eqref{eeq1:th}.
\end{proof}

\begin{proof}[\textbf{Proof of Lemma \ref{lem:low}}]
We begin by noticing that 
$$ \inf_{\hat T\in \mathbb{R}^p}\mathbb{E}_{\pi}\mathbf{E}_\beta \sum_{j=1}^p | {\hat T}_j(X,Y) - \beta_j |^{q}  = \sum_{j=1}^{p}\inf_{\hat T_{j}\in \mathbb{R}}\mathbb{E}_{\pi}\mathbf{E}_\beta  | {\hat T}_j(X,Y) - \beta_j |^{q}.  $$
It is easy to check that
\begin{equation}\label{eq:inv}
    \forall a \in \mathbb{R}^{p}, \forall j=1,\dots,p \quad \inf_{\hat T_{j}\in \mathbb{R}}\mathbb{E}_{\pi}\mathbf{E}_\beta  | {\hat T}_{j}(X,Y) - \beta_j |^{q}  = \inf_{\hat T_{j}\in \mathbb{R}}\mathbb{E}_{\pi}\mathbf{E}_\beta  | {\hat T}_{j}(X,Y-a) - \beta_j |^{q}.
\end{equation}
Using conditioning, one may also notice that 
\begin{equation}\label{eq:cond}
    \inf_{\hat T_{j}\in \mathbb{R}}\mathbb{E}_{\pi}\mathbf{E}_\beta  | {\hat T}_j(X,Y) - \beta_j |^{q} \geq \mathbb{E}_{\pi_{-j}}  \left( \inf_{\hat T_{j}\in \mathbb{R}}\mathbb{E}_{\pi_{j}}\mathbf{E}_\beta | {\hat T}_j(X,Y) - \beta_j |^{q} \Big| \beta_{-j}\right),
\end{equation}
where $\beta_{-j}$ represents the vector $\beta$ deprived of $\beta_{j}$ and $\pi_{-j}$ the corresponding prior.
Hence, we get from (\ref{eq:inv}) and (\ref{eq:cond}) that
$$ \inf_{\hat T_{j}\in \mathbb{R}}\mathbb{E}_{\pi}\mathbf{E}_\beta  | {\hat T}_j(X,Y) - \beta_j |^{q} \geq \mathbb{E}_{\pi_{-j}}  \left( \inf_{\hat T_{j}\in \mathbb{R}}\mathbb{E}_{\pi_{j}}\mathbf{E}_\beta | {\hat T}_j(X,\tilde{Y}_{j}) - \beta_j |^{q} \Big| \beta_{-j}\right), $$
where $\tilde{Y}_{j} = Y - \sum_{i\neq j} X_{i}\beta_{i} = \beta_{j}X_{j} + \sigma \xi$. We remove the last conditional expectation and replace the dependence on $X$ by $X_{j}$, since the observable $\tilde{Y}_{j}$ depends only on $\beta_{j}$ and $X_{j}$.
\end{proof}

\begin{proof}[\textbf{Proof of Theorem \ref{th:low2}}]
We apply Theorem \ref{th:low1} with $\Theta_{s,a} = \Omega_{s,a}$. Let $S$ a support of size $ s$, and consider the prior $\beta$ such that $\beta_{S^{c}} = 0 $ and $\beta_{S} = Z$,
where $Z \in \mathbb{R}^{s}$ is a Gaussian random vector distributed following $\mathcal{N}(\mu,\nu^{2}\mathbb{I}_{s})$ where $\mu, \nu >0$ are defined later. We have
$$\inf_{\hat{\beta}} \sup_{ \beta \in \beta_{s,a}} \mathbf{E}_\beta |\hat{\beta} - \beta|^{q}
\geq  \inf_{\hat T\in \mathbb{R}^p}\mathbb{E}_{\pi}\mathbf{E}_\beta \sum_{j=1}^p | {\hat T}_j(X,Y) - \beta_j |^{q} - C_{q}\,\mathbb{E}_{\pi}\left[ \left( \mathbf{E}(\|\beta^{A}\|^{q}_{q}|Y) + \|\beta\|_{q}^{q} \right)\mathbf{1}( \beta \not\in \Omega_{s,a} ) \right].$$
We first upper-bound the second term
$$ \mathbb{E}_{\pi}\left[ \left( \mathbf{E}^{A}(\|\beta^{A}\|^{q}_{q}|Y) + \|\beta\|_{q}^{q} \right)\mathbf{1}( \beta \not\in \Theta_{s,a} ) \right] \leq 2\sqrt{\mathbb{E}_{\pi} \|\beta\|^{2q}_{q}}\sqrt{\mathbf{P}( \beta \not\in \Theta{s,a} )},   $$
since $\|\beta^{A}\|^{q}_{q} \leq \|\beta\|^{q}_{q}$. It is easy to check that for some $C > 0$ we have
$$  \mathbf{P}( \beta \not\in \Theta_{s,a} ) \leq s\mathbf{P}( |\beta_{1}| \leq a ) \leq C s e^{-\frac{(\mu_{1} - a)^{2}_{+}}{2\nu^{2}}}.$$
By choosing $\mu_{1} = a + \nu^{2} $, we get for some $C_{q}>0$
$$ \mathbb{E}_{\pi}\left[ \left( \mathbb{E}^{A}(\|\beta^{A}\|^{q}_{q}|Y) + \|\beta\|_{q}^{q} \right)\mathbf{1}( \beta \not\in \Theta_{s,a} ) \right] \leq C_{q}\sqrt{s}p\sqrt{a^{2q} + \nu^{4q}+\nu^{2q}}e^{-\frac{\nu^{2}}{2}}. $$
Using lemma \ref{lem:low} combined with Anderson lemma for Gaussian priors we get
$$ \inf_{\hat T_{j}\in \mathbb{R}}\mathbb{E}_{\pi_{j}}\mathbf{E}_\beta | {\hat T}_j(X,\tilde{Y}_{j}) - \beta_j |^{q} = \mathbf{E} \left( \left(\frac{\nu \sigma}{\nu \|X_{j}\|_{2} + \sigma} \right)^{q}|\xi_{1}|^{q} \right). $$
We conclude that $\forall \nu>0$, we have

$$\underset{\hat{\beta}}{\inf}\underset{\beta \in \Omega^{p}_{s,a}}{\sup} \mathbf{E}_{\beta}\left( \| \hat{\beta} - \beta \|_{q}^{q}\right) \geq \sum_{j\in S}\left(\frac{\nu \sigma}{\nu \|X_{j}\|_{2} + \sigma} \right)^{q}\mathbf{E}\left(|\xi_{1}|^{q} \right)- C_{q}\sqrt{s}p\sqrt{a^{2q} + \nu^{4q}+\nu^{2q}}e^{-\frac{\nu^{2}}{2}}.$$

The result follows by taking the limit $\nu \to \infty$.
\end{proof}

\begin{proof}[\textbf{Proof of Theorem \ref{th:low3}}]
We are going to mimic the previous proof using a different prior. We apply Theorem \ref{th:low1} with $\Theta_{s,a} = \Omega_{s,a}$. Consider the prior $\beta$ such that $\beta = a \eta$,
where $\eta \in \{0,1\}^{p}$ be a Bernoulli random vector with i.i.d entries and $\mathbf{E}(\eta_{i}) = \frac{s'}{p}$, $s'\in(0,s)$. We have
$$\inf_{\hat{\beta}} \sup_{ \beta \in \Theta{s,a}} \mathbf{E}_\beta |\hat{\beta} - \beta|^{q}
\geq  \inf_{\hat T\in \mathbb{R}^d}\mathbb{E}_{\pi}\mathbf{E}_\beta \sum_{j=1}^p | {\hat T}_j(X,Y) - \beta_j |^{q} - C_{q}\,\mathbb{E}_{\pi}\left[ \left( \mathbf{E}(\|\beta^{A}\|^{q}_{q}|Y) + \|\beta\|_{q}^{q} \right)\mathbf{1}( \beta \not\in \Theta{s,a} ) \right].$$
First notice that in this case
$$   \beta \in \Theta{s,a} \quad \text{if and only if} \quad  |\eta|_{0} \leq s. $$
Hence $\|\beta^{A}\|^{q}_{q} \leq a^{q}|\eta|_{0} \leq sa^{q}$.
We first upper-bound the second term
$$ \mathbb{E}_{\pi}\left[ \left( \mathbf{E}^{A}(\|\beta^{A}\|^{q}_{q}|Y) + \|\beta\|_{q}^{q} \right)\mathbf{1}( \beta \not\in \Theta_{s,a} ) \right] \leq a^{q}\mathbb{E}_{\pi}\left[ 2 |\eta|_{0} \mathbf{1} ( |\eta|_{0} \geq s+1 ) \right], $$
since $|\eta|_{0} > s$.
Using same arguments as in \cite{butucea2018}, we conclude that
$$ \mathbb{E}_{\pi}\left[ \left( \mathbf{E}^{A}(\|\beta^{A}\|^{q}_{q}|Y) + \|\beta\|_{q}^{q} \right)\mathbf{1}( \beta \not\in \Theta_{s,a} ) \right] \leq 2a^{q}s'e^{-\frac{(s-s')^{2}}{2s}}. $$
Going back to the first term, we get the following lower bound using Lemma \ref{lem:low}
$$ \inf_{\hat T_{j}\in \mathbb{R}}\mathbb{E}_{\pi_{j}}\mathbf{E}_\beta | {\hat T}_j(X,\tilde{Y}_{j}) - \beta_j |^{q} = a^{q}\inf_{\hat T_{j}\in \mathbb{R}}\left( \frac{s'}{p}\mathbf{E}_a | {\hat T}_j(X,\tilde{Y}_{j}) - 1 |^{q} + (1-\frac{s'}{p})\mathbf{E}_0 | {\hat T}_j(X,\tilde{Y}_{j}) |^{q} \right)$$
Minimizing the posterior risk, the Bayes rule gives
$$ \forall q >1,\quad  T_{j}^{*}(X,\tilde{Y}_{j}) =  \frac{1}{1+e^{\frac{a}{q-1}(t_{j}(a)\|X_{j}\|_{2}^{2}-\langle\tilde{Y}_{j},X_{j} \rangle)}}, $$
and for $q=1$ we get
$$  T_{j}^{*}(X,\tilde{Y}_{j}) = \mathbf{1}( \langle\tilde{Y}_{j},X_{j} \rangle \geq t_{j}(a)\|X_{j}\|^{2}_{2} ). $$
Hence we deduce that
$$ \inf_{\hat T_{j}\in \mathbb{R}}\mathbb{E}_{\pi_{j}}\mathbf{E}_\beta | {\hat T}_j(X,\tilde{Y}_{j}) - \beta_j |^{q} \geq \frac{a^{q}}{2^{q}}\left( \frac{s'}{p} \mathbf{P}_{a}( \langle\tilde{Y}_{j},X_{j} \rangle \leq t_{j}(a)\|X_{j}\|^{2} ) + (1-\frac{s'}{p}) \mathbf{P}_{0}( \langle\tilde{Y}_{j},X_{j} \rangle \geq t_{j}(a)\|X_{j}\|_{2}^{2}) \right).$$
Notice that for $q=1$ the term $2^{q}$ is not needed. Replacing $\tilde{Y}_{j}$ by its expression, we recover the lower bound
$$ \Psi(p,s',a,\sigma,X) = \sum_{j=1}^{p}\left(\frac{s'}{p}\mathbf{P}(\varepsilon \geq (a-t_{j}(a))\|X_{j}\|_{2})  + (1-\frac{s'}{p})\mathbf{P}(\varepsilon \geq t_{j}(a)\|X_{j}\|_{2})  \right). $$
Following the proof of \cite{ndaoud2018}, we may use the fact that $s \to \frac{\Psi(s)}{s}$ is decreasing to conclude the proof.
\end{proof}

\begin{proof}[\textbf{Proof of Proposition \ref{prop:low}}]
Combining Theorem \ref{th:low2} and Theorem \ref{th:low3} with $s'=s/2$, we get 
$$ \underset{\hat{\beta}}{\inf}\underset{\beta \in \Omega^{p}_{s,a}}{\sup} \mathbf{E}_{\beta}\left( \| \hat{\beta} - \beta \|_{q}^{q}\right) \geq s\sigma^{q}_{q} \vee \left( \frac{a^{q}}{2^{q+1}} \psi(p,s,a,\sigma)- sa^{q}e^{-s/8} \right). $$
We remind the reader the notation $t^{*} := \sqrt{2\sigma^{2}\log\left(\frac{p}{s}-1\right)}$. In order to prove the result, we handle several cases.
 \begin{itemize}
     \item case $a \geq 10 t^{*}$:\\
     It is easy to check that $a - t(a) \geq a / 4$ and that $t(a) \geq a / 4 + t^{*}$. Hence
      $$ a^{q}\psi(p,s,a,\sigma) \leq Csa^{q}e^{-a^{2}/32\sigma^{2}} \leq C_{q}s.$$
      This shows that the term $\sigma^{q}_{q}s$ is dominating. As a result
      $$ s\sigma^{q}_{q} \vee \left( \frac{a^{q}}{2^{q+1}} \psi(p,s,a,\sigma)- sa^{q}e^{-s/8} \right) \asymp s, $$
      and
      $$ s\sigma^{q}_{q} \vee  a^{q} \psi(p,s,a,\sigma) \asymp s.$$
      This suffises to prove the lower bound.
     \item case $t^{*} \leq a \leq 10 t^{*}$:\\
     Since $s \geq 8q\log\log{p}$, then
     $$ a^{q}e^{-s/8} \leq C_{q}a^{-q/2} \leq C_{q'}. $$
     This leads to 
     $$ \underset{\hat{\beta}}{\inf}\underset{\beta \in \Omega^{p}_{s,a}}{\sup} \mathbf{E}_{\beta}\left( \| \hat{\beta} - \beta \|_{q}^{q}\right) \geq s\sigma^{q}_{q} \vee \left( \frac{a^{q}}{2^{q+1}} \psi(p,s,a,\sigma)- sC_{q'} \right). $$
     We conclude by noticing that $a\vee b \asymp a \vee (b-a)$ for $a,b \geq 0$.
     \item case $a \leq t^{*}$:\\
We observe that $t(t^{*})=t^{*}$. In this case 
$$ \underset{\hat{\beta}}{\inf}\underset{\beta \in \Omega^{p}_{s,a}}{\sup} \mathbf{E}_{\beta}\left( \| \hat{\beta} - \beta \|_{q}^{q}\right) \geq \underset{\hat{\beta}}{\inf}\underset{\beta \in \Omega^{p}_{s,t^{*}}}{\sup} \mathbf{E}_{\beta}\left( \| \hat{\beta} - \beta \|_{q}^{q}\right) \geq \frac{1}{2^{q+1}} st^{*q}\mathbf{P}\left( \sigma \varepsilon \geq 0 \right) - C_{q'}s \geq C_{q''}st^{*q}.$$
Hence
$$ \underset{\hat{\beta}}{\inf}\underset{\beta \in \Omega^{p}_{s,a}}{\sup} \mathbf{E}_{\beta}\left( \| \hat{\beta} - \beta \|_{q}^{q}\right) \geq C_{q''} \sigma^{q}s\log\left(\frac{p}{s}-1\right)^{q/2}.$$
 \end{itemize}
\end{proof}
\begin{proof}[\textbf{Proof of Theorem \ref{thm:up1}}]
Let $\beta$ be a vector in $\Omega_{s,a}^{p}$, we have
$$ \| \beta^{a} - \beta \|_{q}^{q} = \sum_{i \in S} \mathbf{E}\left| \hat{\beta}_{i}^{a}-\beta_{i} \right|^{q} + \sum_{i \in S^{c} }\mathbf{E}\left| \hat{\beta}_{i}^{a}-\beta_{i} \right|^{q}. $$
On $S^{c}$, we have
$$ \hat{ \beta }^{a}_{i} - \beta_{i} = \xi_{i} \mathbf{1}_{ \{|\xi_{i}|> t(a \vee t^{*})\} }.  $$
Hence we get that
$$\mathbf{E}\left| \hat{\beta}_{i}^{a}-\beta_{i} \right|^{q} = \mathbf{E}\left( |\xi_{i}|^{q} \mathbf{1}_{\{ |\xi_{i}|>t(a \vee t^{*}) \} }\right).
$$
Using integration by parts and induction we get  
$$ \forall q \geq 0,\quad \mathbf{E}\left( |\xi_{i}|^{q} \mathbf{1}_{\{ |\xi_{i}|>t(a \vee t^{*}) \} }\right) \leq C_{q}(t(a \vee t^{*})^{q} + \sigma^{q})\mathbf{P}\left( |\xi_{i}| \geq t(a \vee t^{*}) \right),$$
where $C_{q}$ is a universal constant depending only in $q$.
Applying this we get
\begin{align*}
            \mathbf{E}\left| \hat{\beta}_{i}^{a}-\beta_{i} \right|^{q} &= \mathbf{E}\left( |\xi_{i}|^{q} \mathbf{1}_{\{ |\xi_{i}|>t(a \vee t^{*}) \} }\right)\\
            &\leq C_{q}(t(a \vee t^{*})^{q} + \sigma^{q})\mathbf{P}\left( |\xi_{i}| \geq t(a \vee t^{*}) \right).
\end{align*}
Hence 
$$  \mathbf{E}\sum_{i \in S^{c}}\left| \hat{\beta_{i}}-\beta_{i} \right|^{q} \leq  2C_{q}|S^{c}|t(a \vee t^{*})^{q}\mathbf{P}\left( \sigma \varepsilon \geq t\left(a \vee t^{*}\right) \right).$$
The last inequality holds since $t(a \vee t^{*}) \geq c \sigma$ for $s \leq p/4$.

On $S$, we have
$$ \hat{ \beta }^{a}_{i} - \beta_{i} = Y_{i} \mathbf{1}_{ \{|Y_{i}|>t(a)\} } - \beta_{i} = -\xi_{i} - Y_{i} \mathbf{1}_{ \{|Y_{i}|\leq t(a)\} }.  $$
Hence and since $|x+y|^{q} \leq 2^{q-1}(|x|^{q}+|y|^{q})$ we get
\begin{align*}
            \mathbf{E}\left| \hat{\beta^{a}_{i}}-\beta_{i} \right|^{q} &\leq 2^{q-1} \sigma_{q}^{q} + 2^{q-1}  \mathbf{E}\left( |Y_{i}|^{q}\mathbf{1}_{\{|Y_{i}| \leq t(a \vee t^{*})\}} \right)\\
            & \leq 2^{q-1} \sigma_{q}^{q} + 2^{q-1} t(a \vee t^{*})^{q} \mathbf{P}\left( |Y_{i}|\leq t(a \vee t^{*}) \right)\\
            & \leq 2^{q-1} \sigma_{q}^{q} + 2^{q-1}t(a \vee t^{*})^{q}\mathbf{P}\left( |\xi_{i}|\geq (a-t(a \vee t^{*}))_{+} \right).
\end{align*}
We get that on $S$ we have
\begin{equation}\label{eq:thm2}
     \mathbf{E}\sum_{i \in S} | \hat{\beta}_{i}^{a}-\beta_{i} |^{q} \leq C_{q} \left( s \sigma_{q}^{q} + t(a)^{q} |S| \mathbf{P}\left( \sigma \varepsilon > (a - t(a \vee t^{*})_{+})\right) \right).
\end{equation}
Since $(a - t(a \vee t^{*})_{+}) \leq (a\vee t^{*} - t(a \vee t^{*})_{+})$, we get 
$$ \mathbf{E} \sum_{i \in S} | \hat{\beta}_{i}^{a}-\beta_{i} |^{q} \leq C_{q}s\sigma_{q}^{q} + C_{q}t(a \vee t^{*})^{q}|S| \mathbf{P}\left( \sigma \varepsilon > (a\vee t^{*} - t(a \vee t^{*})_{+})\right). $$
We conclude that
$$ \mathbf{E}\left( \|\hat{\beta}^{a} - \beta \|_{q}^{q} \right) \leq C_{q}\sigma_{q}^{q}s  +  C_{q}t(a\vee t^{*})^{q}\psi_{+}(p,s,t^{*} \vee a,\sigma).  $$
Hence for $a \geq t^{*}$, the result is immediate, since $t(a \vee t^{*}) \leq t(a) \leq a$. For $a < t^{*}$ we have 
$$ \mathbf{E}\left( \|\hat{\beta}^{a} - \beta \|_{q}^{q} \right) \leq C_{q}\sigma_{q}^{q}s  +  C_{q}\sigma^{q}\log( \frac{p}{s}-1)^{q/2} \psi_{+}(p,s,t^{*},\sigma). $$
It is easy to verify 
$$ \psi_{+}(p,s,t^{*},\sigma) \leq s + (p-s)  \frac{s}{p-s} \leq 2s, $$
and hence 
$$ \mathbf{E}\left( \|\hat{\beta}^{a} - \beta \|_{q}^{q} \right) \leq C_{q} \left( \sigma^{q}s\log(\frac{p}{s}-1)^{q/2} + \sigma_{q}^{q}s\right) \leq C_{q'}\sigma^{q}s\log(\frac{p}{s}-1)^{q/2}, $$
since $s \leq \frac{p}{4}$ and $\log(\frac{p}{s}-1) \geq 1$.
\end{proof}

\begin{proof}[\textbf{Proof of Theorem \ref{thm:up2}}]
Let us first notice that for $\epsilon \in [0,1]$ we have
$$ t(a_{q}(\epsilon)) = \sqrt{2\sigma^{2}\log(\frac{p}{s}-1) + q \epsilon \sigma^{2}\log\log(\frac{p}{s}-1)}, $$
and 
$$ a_{q}(\epsilon) - t(a_{q}(\epsilon)) = \sqrt{ q \epsilon \sigma^{2}\log\log(\frac{p}{s}-1)}.  $$
Following the previous proof we have
$$  \mathbb{E}\sum_{i \in S^{c}}\left| \hat{\beta_{i}}-\beta_{i} \right|^{q} \leq  2C_{q}pt(a)^{q}\mathbf{P}\left( \sigma \varepsilon > t\left(a\right) \right).$$
Since $\epsilon \in [0,1]$ we have that
$ t(a_{q}(\epsilon)) \leq \sqrt{2\sigma^{2}\log(\frac{p}{s}-1)}(1+o(1)).$
Moreover 
$$  \mathbf{P}\left( \sigma \varepsilon > t\left(a\right) \right) \leq C\sigma\frac{e^{-t(a)^{2}/2\sigma^{2}}}{t(a)} \\
\leq  \frac{ C\sigma }{t(a)} \frac{s}{p-s} \frac{1}{\log(p/s-1)^{q\epsilon/2}}.$$
Hence
$$  \mathbf{E}\sum_{i \in S^{c}}\left| \hat{\beta_{i}}-\beta_{i} \right|^{q} \leq  C_{q}\frac{s\log(p/s-1)^{\frac{q}{2}(1-\epsilon)}}{\sqrt{\log(p/s-1)}}.$$
We can now notice that on $S^{c}$ we have
$$  \mathbf{E}\sum_{i \in S^{c}}\left| \hat{\beta_{i}}-\beta_{i} \right|^{q} = \underset{\frac{s}{p} \to 0}{o}\left( \Phi_{o}\right).$$
In order to prove the Theorem we focus on the error in the support. Remember that on $S$ we have
$$ \hat{ \beta }^{a}_{i} - \beta_{i} = Y_{i} \mathbf{1}_{ \{|Y_{i}|>t(a)\} } - \beta_{i} = -\xi_{i} - Y_{i} \mathbf{1}_{ \{|Y_{i}|\leq t(a)\} }.  $$
\begin{itemize}
    \item case $a \leq a(0)$:\\
    In this case $a(0) = t(a(0)) = \sqrt{2\sigma^{2}\log(p/s-1)}$.
    We use the following inequality 
    $$ \forall a,b\in \mathbb{R}, q\geq 1, \quad |a+b|^{q} \leq |a|^{q} + q|a+b|^{q-1}|b|. $$
    Hence 
    \begin{align*}
        | \xi_{i} - Y_{i} \mathbf{1}_{ \{|Y_{i}|\leq t(a)\} }|^{q} &\leq |Y_{i} \mathbf{1}_{ \{|Y_{i}|\leq t(a)\} }|^{q} + q| \xi_{i}|| \xi_{i} - Y_{i} \mathbf{1}_{ \{|Y_{i}|\leq t(a)\} }|^{q-1} \\
    &\leq |Y_{i} \mathrm{1}_{ \{|Y_{i}|\leq t(a)\} }|^{q} + q| \xi_{i}|2^{q}\left(| \xi_{i}|^{q-1} + | Y_{i} \mathbf{1}_{ \{|Y_{i}|\leq t(a)\} }|^{q-1} \right) \\
    & \leq t(a)^{q} + q2^{q}\left(| \xi{i}|^{q} + |\xi_{i}|t(a)^{q-1} \right).
    \end{align*}
    As a consequence 
    $$ \mathbf{E}|\hat{ \beta }^{a}_{i} - \beta_{i}|^{q} \leq t(a)^{q} + q2^{q}(\sigma_{q}^{q} + \sigma_{1}t(a)^{q-1} )\leq t(a)^{q}(1+o(1)). $$
    The last inequality holds since $t(a) \to \infty $ as $s/p \to 0$. We conclude that
    $$ \sum_{i\in S} \mathbf{E}|\hat{ \beta }^{a}_{i} - \beta_{i}|^{q} \leq \Phi_{o}(1+o(1)).$$
    \item case $a = a(\epsilon)$ for $\epsilon \in (0,1)$:\\
    In this case and following same steps in previous case
    $$ | \xi_{i} - Y_{i} \mathrm{1}_{ \{|Y_{i}|\leq t(a)\} }|^{q} \leq  t(a)^{q}\mathbf{1}_{ \{|Y_{i}|\leq t(a)\} } + q2^{q}\left(| \xi_{i}|^{q} + |\xi_{i}|t(a)^{q-1}\mathbf{1}_{ \{|Y_{i}|\leq t(a)\} } \right).$$
    Remember that
$$ \forall q \geq 0,\quad \mathbf{E}\left( |\xi_{i}|^{q} \mathrm{1}_{\{ |\xi_{i}|>t(a) \} }\right) \leq C_{q}(t(a)^{q} + \sigma_{2}^{q})\mathbf{P}\left( |\xi_{i}| \geq t(a) \right).$$
Hence 
$$ \mathbf{E}(|\xi_{i}|\mathbf{1}_{ \{|Y_{i}|\leq t(a)\} } ) \leq  \mathbf{E}(|\xi_{i}|\mathrm{1}_{ \{|\xi_{i}|\geq a - t(a)\} } ) \leq ((a-t) + \sigma)\mathbf{P}(|\xi_{i}| \geq a-t)\leq \log(t)\mathbf{P}(|Y_{i}|\leq t).$$
  
    We get that
    $$ \mathbf{E}|\hat{ \beta }^{a}_{i} - \beta_{i}|^{q} \leq t^{*q}(1+o(1))\mathbf{P}(|Y_{i}| \leq t(a)) + C_{q}\sigma_{q}^{q}. $$
    One may notice that
    $$ \mathbf{P}(|Y_{i}| \leq t(a)) \leq 2\mathbf{P}(\sigma \epsilon \geq (a-t(a))_{+}). $$ 
    Using the Gaussian tail inequality and the fact that $\epsilon>0$, we get 
    $$ \mathbf{P}(|Y_{i}| \leq t(a)) \leq \frac{t^{*-q\epsilon}}{1+\sqrt{\frac{\pi}{2}}q\epsilon\log\log(p/s-1)}(1+o(1)). $$
    Since $t^{*q(1-\epsilon)}/\log(t) \to \infty$ we conclude that
    $$ \sum_{i\in S} \mathbf{E}|\hat{ \beta }^{a}_{i} - \beta_{i}|^{q} \leq \frac{st^{*q(1-\epsilon)}}{1+\sqrt{\frac{\pi}{2}}q\epsilon\log\log(p/s-1)}(1+o(1)).$$
     \item case $a \geq a(1)$:\\
     In this case it suffices to prove the result for $a=a(1)$ since the minimax risk in increasing with respect to $a$.
\begin{align*}
    | \xi_{i} - Y_{i} \mathbf{1}_{ \{|Y_{i}|\leq t(a)\} }|^{q} &\leq |\xi_{i} |^{q} + q| Y_{i} 
    \mathbf{1}_{ \{|Y_{i}|\leq t(a)\} } || \xi_{i} - Y_{i} \mathbf{1}_{ \{|Y_{i}|\leq t(a)\} }|^{q-1} \\
    &\leq |\xi_{i} |^{q} + C_{q}\left( t^{*}|\xi_{i}|^{q-1}\mathbf{1}_{ \{|Y_{i}|\leq t(a)\} } + t^{*q}\mathbf{1}_{ \{|Y_{i}|\leq t(a)\} }   \right).
\end{align*}  
In the previous case we proved that
$$ \mathbf{E}(|\xi_{i}|^{q-1}\mathbf{1}_{ \{|Y_{i}|\leq t(a)\} } ) \leq \log(t^{*})^{q-1}\mathbb{P}(|Y_{i}|\leq t(a)),$$
and that
$$ \mathbf{P}(|Y_{i}| \leq t(a)) \leq C\frac{t^{*-q}}{\log\log(p/s)+1}. $$
Hence 
$$ \mathbf{E}\left( t^{*}|\xi_{i}|^{q-1}\mathbf{1}_{ \{|Y_{i}|\leq t(a)\} } + t^{*q}\mathbf{1}_{ \{|Y_{i}|\leq t(a)\} }   \right) \leq \frac{C}{\log\log(p/s)}=o(\sigma_{q}^{q}). $$
It follows that
$$ \sum_{i\in S} \mathbf{E}|\hat{ \beta }^{a}_{i} - \beta_{i}|^{q} \leq s\sigma^{q}_{q}(1+o(1)).$$
This concludes the proof of this theorem.
\end{itemize}
\end{proof}

\begin{proof}[\textbf{Proof of Corollary \ref{cor:last}}]
Based on the fact that
$$
\underset{\hat{\beta}}{\inf}\underset{\beta \in \Omega^{p}_{s,a}}{\sup}  \mathbf{E}_{\beta}\left( \| \hat{\beta} - \beta \|_{q}^{q}\right) \geq s\sigma_{q}^{q},$$
we observe that the second result is a direct consequence of Theorem \ref{thm:up2}.
In order to conclude, we need to show that for $a < a_{q}(1)$ we have
$$
\underset{\hat{\beta}}{\inf}\underset{\beta \in \Omega^{p}_{s,a}}{\sup} \frac{ \mathbf{E}_{\beta}\left( \| \hat{\beta} - \beta \|_{q}^{q}\right)}{ \Phi_{o}(a) } \geq \frac{1}{4}+o(1).
$$
In what follows, we assume that $a < a_{q}(1)$. Going back to the initial lower bound with $s'=s/2$, and using the fact that $s \geq 8q\log\log{p}$, we have
$$ \quad \underset{\hat{\beta}}{\inf}\underset{\beta \in \Omega^{p}_{s,a}}{\sup} \mathbf{E}_{\beta}\left( \| \hat{\beta} - \beta \|_{q}^{q}\right) \geq\frac{1}{2} sa^{q}\mathbf{E}\left( T_{1}^{q} \right) - C_{q''}sa(1)^{q}a(0)^{-2q},$$
where
$$ \forall q >1,\quad  T_{1} =  \frac{1}{1+e^{-\frac{a}{q-1}(t(a) - a + \xi_{1})}},$$
and for $q=1$
$$
T_{1} = \mathbf{1}(\xi_{1} \geq a - t(a)).
$$
Since
$$ sa(1)^{q}a(0)^{-2q} = o(s), $$
we get immediately that
$$ sa(1)^{q}a(0)^{-2q} = o(\Phi_{o}(a)).$$
It is sufficient to prove that
$$ sa^{q}\mathbf{E}\left( T_{1}^{q} \right) \geq \frac{1}{2}\Phi_{o}(a)(1+o(1)). $$
For $q=1$, we have $\mathbf{E}(T_{1}) = \mathbf{P}(\xi_{1} \geq a - t(a)) \geq \mathbf{P}(|\xi_{1}| \geq (a - t(a))_{+}) / 2$. 
Using the fact that the Gaussian tail bounds presented in Appendix are sharp combined with the proof of the previous upper bound we can verify that for $q=1$
$$ sa\mathbf{E}\left( T_{1} \right) \geq \frac{1}{2}\Phi_{o}(a)(1+o(1)). $$
For $q>1$ it is enough to prove that
$$ \mathbf{E}(T_{1}^{q}) \geq \mathbf{P}(\xi_{1} \geq a - t(a))(1+o(1)). $$
For $a = a_{q}(0)$ we have that $t(a) = a$, hence
$$
\mathbf{E}(T_{1}^{q}) 
= \mathbf{E}\left( \left(\frac{1}{1+e^{\frac{a}{q-1} \xi_{1}}} \right)^{q} \right) \to \mathbf{P}(\xi_{1} \geq 0).$$
The limit is a consequence of the dominated convergence theorem and proves the result.
The last case is when $a = a_{q}(\epsilon)$ with $\epsilon \in (0,1)$. Let us just recall that $a \asymp \sqrt{\log(p/s)}$ and $a-t(a) \asymp \sqrt{\log\log(p/s)}$. Let $\alpha_{s}>0$ be a sequence satisfying $$\alpha_{s}.a \to \infty \text{ and } \alpha_{s}.(a-t) \to 0. $$
We have 
$$ \mathbf{E}(T_{1}^{q}) \geq \mathbf{E}\left(\left( \frac{1}{1+e^{\frac{-a\alpha_{s}}{q-1} }} \right)^{q} \mathbf{1}(\xi_{1} \geq a - t(a) + \alpha_{s}) \right) \geq \left(\frac{1}{1+e^{\frac{-a\alpha_{s}}{q-1} }}\right)^{q} \mathbf{P}(\xi_{1} \geq a - t(a) + \alpha_{s}).  $$
Using the monotony of cumulative distribution functions, we get that
$$
\mathbf{E}(T_{1}^{q}) \geq (1+o(1))\left(\mathbf{P}(\xi_{1} \geq a - t(a)) - Ce^{-(a-t(a)+\alpha_{s})^{2}/2\sigma^{2}}\alpha_{s} \right). 
$$
Using the limiting behaviour of $\alpha_{s}$ we get 
$$ \mathbf{E}(T_{1}^{q}) \geq \mathbf{P}(\xi_{1} \geq a - t(a))(1+o(1)). $$
This concludes the proof.
\end{proof}
\begin{proof}[\textbf{Proof of Theorem \ref{thm:up1ad}}]
First notice that $t_{s}^{*} = t(a(1))$ and that $\hat{\beta}^{s} = \hat{\beta}^{a(1)}$. Hence and using Theorem \ref{thm:up2} we get for all $a \geq a(1)$  
$$  \underset{\beta \in \Omega^{p}_{s,a}}{\sup} \mathbf{E}_{\beta}\left( \| \hat{\beta}^{s} - \beta \|_{q}^{q}\right)\leq \underset{\beta \in \Omega^{p}_{s,a(1)}}{\sup} \mathbf{E}_{\beta}\left( \| \hat{\beta}^{s} - \beta \|_{q}^{q}\right) \leq  C_{q}s\sigma_{q}^{q}.$$
On $S$, we have
$$ \hat{ \beta }^{s}_{i} - \beta_{i} = Y_{i} \mathbf{1}_{ \{|Y_{i}|>t^{*}_{s}\} } - \beta_{i} = -\xi_{i} - Y_{i} \mathbf{1}_{ \{|Y_{i}|\leq t^{*}_{s}\} }.  $$
Hence
$$\mathbf{E}_{\beta}\left( \sum_{i \in S} | \hat{\beta}_{i}^{s} - \beta_{i} |^{q}\right) \leq C_{q}(s\sigma_{q}^{q} + st^{*q}_{s}).
$$
On $S^{c}$, and since $t_{s}>t^{*}$, it is easy to check using previous proofs that
$$  
\mathbf{E}_{\beta}\left( \sum_{i \in S^{c}} | \hat{\beta}_{i}^{s} - \beta_{i} |^{q}\right) \leq
C'_{q} st^{*q}_{s}.
$$
This concludes the proof.
\end{proof}
\vskip 0.2in
\newpage
\bibliographystyle{unsrt}

\end{document}